
\input amstex 
\documentstyle{amsppt}
\input bull-ppt
\keyedby{bull309e/mhm}
\define\IM{\operatorname{Im}}
\define\ess{\operatorname{ess}}
\define\Ac{\operatorname{ac}}
\define\sc{\operatorname{sc}}

\define\N{\Bbb N}

\define\C{\Bbb C}
\define\R{\Bbb R}

\define\PP{\Cal P}
\define\BB{\Cal B}

%

\topmatter
\cvol{27}
\cvolyear{1992}
\cmonth{October}
\cyear{1992}
\cvolno{2}
\cpgs{266-272}
\title New Types of Soliton Solutions \endtitle
\author F. Gesztesy, W. Karwowski, and Z. Zhao \endauthor
\address Department of Mathematics, University of 
Missouri, Columbia, Missouri
65211\endaddress
\ml mathfg\@umcvmb.bitnet \endml
\address Institute of Theoretical Physics, University of 
Wroclaw, 50-205
Wroclaw, Poland\endaddress
\ml iftuwr\@plwrtu11 \endml
\address Department of Mathematics, University of 
Missouri, Columbia, Missouri 
65211\endaddress
\ml mathzz\@umcvmb.bitnet \endml
\date November 19, 1991\enddate
\subjclassrev{Primary 35Q51, 35Q53; Secondary 58F07}
\thanks The first two authors were partially supported by 
the BiBoS-Research
Center at the Faculty of Physics, University of Bielefeld 
and the Department of
Mathematics, University of Bochum, FRG.  The first author 
was also partially
supported by the Norwegian Research Council for Science 
and the Humanities
during a stay at the Division of Mathematical Sciences of 
the University of
Trondheim, Norway.  The second author was also partially 
supported by
CNRS-Marseille, the University of Provence I, Marseille, 
France, and the
German-Polish scientific collaboration program\endthanks
\abstract We announce a detailed investigation of
limits of N-soliton solutions of the Korteweg-deVries 
(KdV) equation as
$N$ tends to infinity.  Our main results provide new 
classes of
KdV-solutions including in particular new types of 
soliton-like
(reflectionless) solutions.  As a byproduct we solve an 
inverse spectral
problem for one-dimensional Schr\"odinger operators and 
explicitly
construct smooth and real-valued potentials that yield a 
purely
absolutely continuous spectrum on the nonnegative real 
axis and give rise
to an eigenvalue spectrum that includes any prescribed 
countable and
bounded subset of the negative real axis.\endabstract
\endtopmatter

\document

\heading Introduction\endheading
In this note we announce the construction of new types of 
soliton-like
solutions of the Korteweg-de Vries (KdV)-equation.  More 
precisely, we
offer a solution to the following problem:

{\it Construct new classes of KdV-solutions by taking 
limits of
$N$-soliton solutions as $N\rightarrow \infty$.}

As it turns out, our solution to this problem is 
intimately connected
with a solution to the following inverse spectral problem 
in connection
with one-dimensional Schr\"odinger operators $H=-\
{d^2}/{dx^2}+V$ in $L^2(\R)$:

{\it Given any bounded and countable subset 
$\{-\kappa^2_j\}_{j\in \N}$ of
$(-\infty,0)$, construct a\/} ({\it smooth and 
real-valued\/}) {\it potential
$V$ such that $H=-{d^2}/{dx^2}+V$ has a purely absolutely 
continuous
spectrum equal to $[0,\infty)$ and the set of eigenvalues 
of $H$ includes the
prescribed set\/} $\{-\kappa^2_j\}_{j\in \N}$.

In addition, we also construct a new class of reflectionless
KdV-solutions in which the underlying Schr\"odinger 
operator has
infinitely many negative eigenvalues accumulating at zero.

Although we present our results exclusively in the 
KdV-context, it will
become clear later on that our methods are not confined to 
KdV-type
equations but are widely applicable in the field of 
integrable systems.

Before formulating our results in detail we briefly review 
some
background material.  The celebrated $N$-soliton solutions 
$V_N(t,x)$ of
the KdV-equation
$$
\text {KdV}(V)=V_t-6VV_x+V_{xxx}=0 \tag 1
$$
described, e.g., in \cite {1, 2, 4}
$$
\gathered 
V_N(t,x)=-2\partial^2_x\ln \det [1_N+C_N(t,x)],\qquad 
(t,x)\in \R^2,\\
C_N(t,x)=\left[\frac
{c_jc_\ell}{\kappa_j+\kappa_\ell}e^{4(\kappa^3_j+
\kappa^3_\ell)
t-(\kappa_j+\kappa_\ell)x}\right]^N_{j,\ell=1},\\
\shoveright { c_j>0,\ \kappa_j>0,\
1\leq j\leq N,\ N\in \N }\endgathered \tag 2
$$
are well known to be isospectral and reflectionless 
potentials $V_N(t,x)$
in connection with the one-dimensional Schr\"odinger 
operator
$$
H_N(t)=-\frac {d^2}{dx^2}+V_N(t,.) \tag 3
$$
in $L^2(\R)$.  In particular, the spectrum
$\sigma(H_N(t))$ of $H_N(t)$ is independent of $t$ and is 
given by
$$
\sigma (H_N(t))=\{-\kappa^2_j\}^N_{j=1}\cup [0,\infty) 
\tag 4
$$
with purely absolutely continuous spectrum $[0,\infty)$.  
Hence $H_N(t)$
are isospectral deformations of $H_N(0)$, which is clear 
from the Lax
formalism connecting (1) and (3).  The reflectionless 
property of $V_N$
manifests itself in the $(t$-independent) scattering 
matrix $S_N(k)$ in
$\C^2$ associated with the pair $(H_N(t), 
H_0)$
$$ 
S_N(k)=\pmatrix T_N(k)&0\\
0&T_N(k)\endpmatrix,\qquad T_N(k)=\prod\limits^N_{j=1}\frac
{k+i\kappa_j}{k-i\kappa_j},\ k\in \C\backslash 
\{i\kappa_j\}^N_{j=1},
$$ 
where $H_0=-{d^2}/{dx^2}$ and $z=k^2$ is the spectral 
parameter
corresponding to $H_0$.  Here $T_N(k)$ denotes the 
transmission
coefficient and the vanishing of the off-diagonal terms in 
$S_N(k)$
exhibits reflection coefficients identical to zero at all 
energies.  In
more intuitive terms this remarkable and highly 
exceptional behavior can
be described as follows:\ \ If one views $V$ as 
representing an
``obstacle'' for an incoming ``signal'' (wave, etc.)\ then 
the outgoing
signal generically consists of two parts, a transmitted 
and a reflected
one.  It is in the exceptional case of $N$-soliton 
potentials $V_N$ such
as (2) that the reflected part of the outgoing signal is 
entirely missing
and hence the obstacle appears to be completely 
transparent independently
of the wavelength of the incoming signal.

Incidentally, (4) offers a solution to the following 
inverse spectral
problem:\ \ Given the finite set 
$\{-\kappa^2_j\}^N_{j=1}\subset
(-\infty,0)$, construct (smooth and real-valued) 
potentials $V_N$ such
that $H_N=-{d^2}/{dx^2}+V_N$ has a purely absolutely 
continuous
spectrum equal to $[0,\infty)$ and precisely the eigenvalues
$\{-\kappa^2_j\}^N_{j=1}$.

A natural generalization of this fact would be to ask 
whether one can
choose a sequence $\{c_j>0\}_{j\in \N}$ such that for an 
arbitrarily
prescribed bounded and countable set 
$\{-\kappa^2_j\}_{j\in \N}\subset
(-\infty,0),\ V_N(t,x)$ converge to a smooth KdV-solution 
$V_\infty(t,x)$
as $N\rightarrow \infty$ with the associated Schr\"odinger 
operator
$H_\infty(t)=  
-{d^2}/{dx^2}+V_\infty(t,.)$ having the purely
absolutely continuous spectrum $[0,\infty)$ and containing 
the set
$\{-\kappa^2_j\}_{j\in \N}$ in its point spectrum.

The main goal of this note is to present an affirmative 
answer to this
question.

\heading Main results\endheading

\thm {Theorem 1} Assume $\{\kappa_j>0\}_{j\in \N}\in
\ell^\infty(\N),\ \kappa_j\neq \kappa_\ell$ for $j\neq 
\ell$, and choose
$\{c_j>0\}_{j\in
\N}$ such that
$\{c^2_j/\kappa_j\}_{j\in \N}\in \ell^1(\N)$.  Then $V_N$ 
converges
pointwise to some $V_\infty\in C^\infty(\R^2)\cap 
L^\infty(\R^2)$ as
$N\rightarrow \infty$ and
\varroster
\item "(i)" \<$\lim_{x\rightarrow +\infty}V_\infty(t,x)=0$ 
and
$$
\lim\limits_{N\rightarrow \infty}\sup\limits_{(t,x)\in K}
|\partial^m_t\partial^n_xV_N(t,x)-\partial^m_t%
\partial^n_xV_\infty(t,x)|=0,
\qquad m,n \in \Bbb N_0, \tag 5
$$
for any compact subset $K\subset \R^2$.  Moreover,
$$
\text {KdV}(V_\infty)=0. \tag 6
$$
\item "(ii)" Denoting $H_\infty(t)=-
{d^2}/{dx^2}+V_\infty(t,.)$ we have
$$
\gather
\sigma_{\ess}(H_\infty(t))=\{-\kappa^2_j\}'_{j\in \N}\cup 
[0,\infty),
\tag 7\\
\sigma_{\Ac}(H_\infty(t))=[0,\infty), \tag 8\\
[\sigma_p(H_\infty(t))\cup \sigma_{\sc}(H_\infty(t))]\cap
(0,\infty)=\emptyset, \tag 9\\
\{-\kappa^2_j\}_{j\in \N}\subseteq 
\sigma_p(H_\infty(t))\subseteq
\overline {\{-\kappa^2_j\}_{j\in \N}}. \tag 10 \endgather
$$\endroster\ethm
\vskip-0.75pc
The spectral multiplicity of $H_\infty(t)$ on $(0,\infty)$
equals two while $\sigma_p(H_\infty(t))$ is simple.  In 
addition, if
$\{\kappa_j\}_{j\in \N}$ is a discrete subset of 
$(0,\infty)$ (i.e., if 0
is its only limit point) then
$$
\gather
\sigma_{\sc}(H_\infty(t))=\emptyset, \tag 11\\
\sigma (H_\infty(t))\cap
(-\infty,0)=\sigma_d(H_\infty(t))=\{-\kappa^2_j\}_{j\in 
\N}. \tag 12
\endgather
$$
More generally, if $\{-\kappa^2_j\}'_{j\in \N}$ is 
countable then (11)
holds.

Here $\overline A$ denotes the closure of $A\subset \R,\ 
A'$ is the
derived set of $A$ (i.e., the set of accumulation points 
of A),
$\N_0=\N\cup \{0\}$, and $\sigma_{\ess}(.),\ 
\sigma_{\Ac}(.),\
\sigma_{\sc}(.),\ \sigma_d(.)$, and $\sigma_p(.)$ denote 
the essential,
absolutely continuous, singularly continuous, discrete, 
and point
spectrum (the set of eigenvalues) respectively.

Under stronger assumptions on $\{\kappa_j\}_{j\in \N}$ we 
obtain

\thm{Theorem 2} Assume $\{\kappa_j\!>\!0\}_{j\in \N}\in
\ell^1(\N),\ \kappa_j\neq \kappa_\ell$ for $j\neq \ell$, and
choose $\{c_j\!>\!0\}_{j\in \N}$ such that 
$\{c^2_j/\kappa_j\}_{j\in
\N}\in
\ell^1(\N)$.  Then in addition to {\rm (5)} and {\rm (6)} 
we have

{\rm (i)}
$$
\lim\limits_{N\rightarrow
\infty}||\partial^m_t\partial^n_xV_N(t,.) 
-\partial^m_t\partial^n_x
V_\infty(t,.)||_p=0,\qquad m,n\in \N_0,\ 1\leq p\leq \infty.
\tag 13
$$

{\rm (ii)}
\vskip-1pc
$$\gather
\sigma_{\ess}(H_\infty(t))=
\sigma_{\Ac}(H_\infty(t))=[0,\infty),\\
\sigma_p(H_\infty(t))\cap
(0,\infty)=\sigma_{\sc}(H_\infty(t))=\emptyset,\\
\sigma_d(H_\infty(t))=\{-\kappa^2_j\}_{j\in \N}\.
\endgather$$

{\rm (iii)} The $(t$-independent\/{\rm )} scattering 
matrix $S_\infty(k)$
in $\C^2$ associated with the pair $(H_\infty(t), H_0)$ is
reflectionless and given by
$$\gather
S_\infty(k)=\pmatrix 
T_\infty(k)&0\\0&T_\infty(k)\endpmatrix,\\
\shoveright 
{T_\infty(k)=\prod\limits^\infty_{j=1}\frac {k+
i\kappa_j}{k-i\kappa_j},\
k\in \C\backslash \{\{i \kappa_j\}_{j\in \N}\cup \{0\}\}.}
\endgather$$
\ethm

While Theorem~1 solves the two problems stated in the 
introduction,
Theorem~2 constructs a new class of reflectionless 
potentials in
connection with one-dimensional Schr\"odinger operators 
involving an
infinite negative point spectrum accumulating at zero.  
Moreover, suppose
$V\in C^\infty(\R^2)$ to be real-valued with 
$\partial^m_xV(t,.)\in
L^1(\R),\ m\in \N$, and either $V(t,.)\in L^1(\R;(1+
|x|^\varepsilon)\,dx)$
for some $\varepsilon>0$ or that $T(k),\ k\in \R\backslash 
\{0\}$, the
transmission coefficient associated with the pair $(H(t)=-
{d^2}/{dx^2}+V(t,.), H_0)$, has an analytic continuation 
into
$\{\C_+\backslash \{i\kappa_j\}_{-\kappa^2_j\in
\sigma_d(H(0))}\}\cup\{k\in
\C\backslash \{0\}\big|\,|k|<\eta,\IM k\leq 0\}$ for some
$\eta>0\ (\C_+=\{k\in \C|\IM k>0\})$.
Introducing the KdV-invariants $\chi_n\in C^\infty(\R^2)$ by
$$
\chi_1=V,\quad \chi_2=-V_x,\quad 
\chi_{n+1}=-\partial_x \, \chi_n
-\sum\limits^{n-1}_{m=1}\chi_{n-m}\chi_m,\qquad n\geq
2,
$$
an extension of the results in \cite {5, 13} yields the 
conservation 
laws (trace relations)
$$\aligned
-\int_{\R} \,dx\chi_{2n+1} (t,x)=& \frac
{2^{2(n+1)}}{2n+1}\sum\limits_{-\kappa^2_j\in
\sigma_d(H)}\kappa^{2n+1}_j\\
&+(-1)^n2^{2(n+1)}\frac 1\pi
\int^\infty_0dk k^{2n}\ln |T(k)|,\qquad n\in \N_0 
\endaligned\tag 14
$$
assuming $\dsize \sum\nolimits_{-\kappa^2_j\in 
\sigma_d(H)}\kappa_j<\infty$
(see \cite {6} for details).  Since $|T(k)|\leq 1 \text { 
for } k>0$ by the
unitarity of the associated scattering matrix, this yields 
the bounds
$$
\align
-\int_{\R} \,dx\chi_{4m+1}(t,x)&\leq \frac
{2^{2(2m+1)}}{4m+1}\sum\limits_{-\kappa^2_j\in
\sigma_d(H)}\kappa_j^{4m+1},\qquad  m\in \N_0, \tag 15\\
-\int_{\R} \,dx\chi_{4m+3}(t,x)&\geq \frac
{2^{2(2m+2)}}{4m+3}\sum\limits_{-\kappa^2_j\in
\sigma_d(H)}\kappa^{4m+3}_j,\qquad  m\in \N_0.\tag 16 
\endalign
$$
For $m=0$ the bound (15) can be found in \cite {12}.  In 
the case where
$\partial^m_xV(t,.) \in L^1(\R;(1\!+\!|x|)\,dx),
m\in \N_0$, and hence
$\sigma_d(H(t))$ is finite, (15) and (16) are discussed, 
e.g., in \cite {8, 10}.
By (14), the bounds (15) and (16) saturate iff 
$|T(k)|=1\text { for
}k>0$, i.e., iff $V$ is reflectionless.  Consequently, the 
bounds (15)
and (16) saturate if $V$ equals the $N$-soliton solutions 
$V_N$ in (2)
and, in particular, if $V$ is an element of our new class of
reflectionless KdV-solutions $V_\infty$ described in 
Theorem~2.

\heading Sketch of proofs \endheading
The hypotheses in Theorem 1 guarantee that $C_N(t,x)$ 
(viewed as an
operator in
$\ell^2(\N))$ converges for any fixed $(t,x)\in \R^2$ in 
trace norm to
some trace class operator $C_\infty(t,x)\in
\BB_1(\ell^2(\N))$ as
$N\rightarrow \infty$ and hence
$$
\lim\limits_{N\rightarrow 
\infty}V_N(t,x)=V_\infty(t,x)=-2\partial^2_x\ln
\det\nolimits_1[1+C_\infty(t,x)],
$$
where $\det\nolimits_1(.)$ denotes the corresponding 
Fredholm
determinant.  A crucial identity in proving (5) and (6) is
$$
V_\infty(t,x)=-4\sum\limits^\infty_{j=1}\kappa_j\psi_{%
\infty,j}(t,x)^2,
\tag 17
$$
where $\{\psi_{\infty,j}(t,x)\}_{j\in \N}$ turn out to be 
the
eigenfunctions of $H_\infty(t)$ corresponding to the 
eigenvalues
$\{-\kappa^2_j\}_{j\in \N}$, determined by
$$
\aligned
\Psi_\infty(t,x)=& [1+
C_\infty(t,x)]^{-1}\Psi^0_\infty(t,x),\\
\Psi^0_\infty(t,x)=&\{c_je^{-\kappa_jx}\}^T_{j\in \N},\\
\Psi_\infty(t,x)=& \{\psi_{\infty,j}(t,x)\}^T_{j\in \N } 
\endaligned \tag 18
$$
in $\ell^2(\N)$.  Equations (17) and (18) are well known 
in the context of $V_N$
and can be obtained by pointwise limits as $N\rightarrow 
\infty$.  Equation (10)
then follows from strong resolvent convergence of $H_N(t)$ 
to
$H_\infty(t)$ and $\sigma_{\ess}(H_\infty(t))\supseteq 
[0,\infty)$ is a
consequence of $V_\infty(t,x)\rightarrow_{x\rightarrow +
\infty}0$.  Next one
constructs the Weyl $m$-functions $m^\pm_\infty(t,z)$ 
associated with
$H^\pm_{\infty,D}(t)$, the restriction of $H_\infty(t)$ to 
the interval
$(0,\pm \infty)$ with a Dirichlet boundary condition at 0. 
 One obtains
$$
\align
m^\pm_\infty(t,z)=&\pm i\sqrt z
\mp [1\mp i\sum\limits^\infty_{j=1}(\sqrt z\pm
i\kappa_j)^{-1}c_j\psi_{\infty,j}(t,0)]^{-1}\\
&\times i\sum\limits^\infty_{j=1}(\sqrt
z \pm i\kappa_j)^{-1}c_j[\partial_x\psi_{\infty,j}(t,0)-%
\kappa_j
\psi_{\infty,j}(t,0)],\qquad z\in \C\backslash \R,
\endalign
$$ 
defining the branch of $\sqrt z$ by 
$\lim_{\varepsilon\downarrow
0}\sqrt {|\lambda|\pm i\varepsilon}=\pm
|\lambda|^{1/2},\ \lim_{\varepsilon\downarrow 0}\sqrt 
{-|\lambda|\pm
i\varepsilon}=i|\lambda|^{1/2}$.  Since 
$m^\pm_\infty(t,.)$ are bounded on
any region of the type $J_{\delta,R_1,R_2}=\{z=\lambda+i\nu
\big|R_1<\lambda<R_2, \ 0<\nu<\delta\},\ 
\delta,R_1,R_2\!>\!0$
and
$$\lim_{\varepsilon\downarrow 0}\big|\IM [m^\pm
_\infty(t,\lambda+i\varepsilon)]\big|$$
is bounded away from zero for
$\lambda\in (R_1,R_2)$, the spectrum of $H_\infty(t)$ in
$(0,\infty)$
is purely absolutely continuous by Theorem~3.1 of \cite 
{11} and hence (9) and
$\sigma_{\Ac}(H_\infty(t))\supset (0,\infty)$ follow.  
Next one proves
the following lemma on the basis of $H^p$-theory, $0<p<1$ 
(see,
e.g., \cite {3}).

\thm{Lemma 3.3} Let $\{a_j\}_{j\in \N},\ \{b_j\}_{j\in
\N}\subset \R,\ \{a_j\}_{j\in \N}\in \ell^1(\N)$.  Then 
there exists a
real-valued function $f$ on $[0,\infty)$ with
$$
m(\{x\in [0,\infty)\big|f(x)=c\})=0 \quad \text {for each 
}c\in \R
$$
such that
$$
\lim\limits_{\varepsilon\downarrow 
0}\sum\limits^\infty_{j=1}\frac
{a_j}{\sqrt {x-i\varepsilon}-b_j}=f(x) \quad \text {for 
}m\text {-a.e. }x\geq 0.
$$
\ethm
(Here $m$ denotes the Lebesgue measure on $\R$.)  
Identifying
$x=-\lambda,\ \lambda<0,\ b_j=\mp \kappa_j,\ 
a_j=c_j\psi_{\infty,j}(t,0)$
resp.\ 
$a_j=c_j[\partial_x\psi_{\infty,j}(t,0)-\kappa_j\psi_{%
\infty,j}(0,t)]$,
Lemma~3 applied to $m^\pm_\infty(t,.)$ yields the 
existence of
real-valued and finite limits of $m^\pm_\infty(t,\lambda+
i\varepsilon)$ for
a.e.\ $\lambda<0$ as $\varepsilon \downarrow 0$.  Thus
$$
\mu^\pm_{\infty,\text {ac}}((-\infty,0))=0,
$$
where $\mu^\pm_{\infty,\Ac}$ is the absolutely continuous 
part (with
respect to $m)$ of the Stieltjes measure generated by the 
spectral
function of $H^\pm_{\infty,D}(t)$.  Consequently
$\sigma_{\Ac}(H^\pm_{\infty,D}(t))$, being the topological 
support of
$\mu^\pm_{\infty, \Ac}$, is contained in $[0,\infty)$.
This together with
$\sigma_{\Ac}(H_\infty(t))=\sigma_{\Ac}(H^+
_{\infty,D}(t))\cup
\sigma_{\Ac}(H^-_{\infty,D}(t))$ yields 
$$\sigma_{\Ac}(H_\infty(t))\cap
(-\infty,0)=\emptyset$$
and hence (8).  The rest of Theorem~1 is
plain.

Finally we turn to Theorem~2.  Due to (17), the fact that 
$\partial_x^2 \,
\psi_{\infty,j}= 
(V_\infty+\kappa^2_j)\psi_{\infty,j}$ and the KdV-equation 
(6)
for $V_\infty$ one can show it suffices to prove (13) for 
$0\leq m\leq 2$ and
$n=0$.  This is accomplished in a series of steps.  First 
one proves the
crucial identity
$$
\int_{\R}\,dx 
V_\infty(t,x)=-4\sum\limits^\infty_{j=1}\kappa_j,
$$
which follows from
\vskip0.5pc
\noindent{\bf Lemma 4.}\  {\it Assume the hypotheses in 
Theorem {\rm 1}.  Then}
$$
\det\nolimits_1[1+C_\infty(0,x)]=1+\sum\limits_{I\in
\PP}a_{\infty,I}e^{-2\sum_{j\in I}\kappa_jx},
$$
{\it where $\PP$ is the family of all finite, nonempty 
subsets of $\N$ and
$a_{\infty,I}>0$ are positive numbers {\rm (}whose precise 
value turns out to
be immaterial for the proof of Theorem {\rm 2)}}\newline
\vskip0.5pc
\noindent and a detailed study of the
asymptotic behavior of $\det\nolimits_1[1+C_\infty(0,x)]$ as
$|x|\rightarrow \infty$.

In the sequel one repeatedly invokes the identity (17) and 
Vitali's
theorem (\cite {9, p.\ 203}).  The rest of Theorem~2 
follows from Theorem~1(ii)
and scattering theory for $L^1(\R)$-potentials.

Detailed proofs can be found in \cite {7}.

We feel that the simplicity of constructing KdV-solutions 
producing these
remarkable spectral (resp.\ scattering) properties 
represents a
significant result that deserves further investigations.  
In particular,
generalizations, replacing the $N$-soliton KdV-solutions 
$V_N$ by $N$-gap
quasi-periodic KdV-solutions and studying the limit 
$N\rightarrow \infty$
involving accumulations of spectral gaps and bands, appear 
to offer a
variety of interesting and challenging problems.

Due to the close resemblance of the determinant structure 
of the
$N$-soliton solutions of (hierarchies of) integrable 
systems such as the
AKNS-class (particularly the nonlinear Schr\"odinger and 
Sine-Gordon
equations), the Toda lattice, and especially the 
Kadomtsev-Petviashvili
equation, the methods in this paper are by no means 
confined to KdV-type
equations but are widely applicable in this field.

\heading Acknowledgments \endheading
We are indebted to N. Kalton and W. Kirsch for several
stimulating discussions.  In particular, we would like to 
thank N. Kalton
for invaluable help in connection with Lemma 3.


\Refs
\rc
\ref\no 1 
\by P. A. Deift
\paper Applications of a commutation formula
\jour Duke Math. J.
\vol 45
\yr 1978
\pages 267--310
\endref

\ref\no 2
\by P. Deift and E. Trubowitz
\paper Inverse scattering on the line
\jour Commun. Pure Appl. Math.
\vol 32
\yr 1979
\pages 121--251
\endref

\ref\no 3
\by P. L. Duren 
\book Theory of $H^p$ Spaces
\publ Academic Press
\publaddr New York
\yr 1970
\endref

\ref\no 4
\by C. S. Gardner, J. M. Greene, M. D. Kruskal, and R. M. Miura
\paper Korteweg-de Vries equation and generalizations, {\rm VI}. Methods for
exact solution
\jour Comm. Pure Appl. Math.
\vol 27
\yr 1974
\pages 97--133
\endref

\ref\no 5
\by I. M. Gel\cprime fand and L. A. Dikii
\paper Asymptotic behavior of the resolvent of Sturm-Liouville equations and the
algebra of the Korteweg-de Vries equations
\jour Russian Math. Surveys
\vol 30:5
\yr 1975
\pages 77--113
\endref

\ref\no 6
\by F. Gesztesy and H. Holden
\paperinfo in preparation
\endref

\ref\no 7
\by F. Gesztesy, W. Karwowski and Z. Zhao
\paper Limits of soliton solutions
\jour Duke Math. J.
\toappear
\endref

\ref\no 8
\by H. Grosse
\paper Quasiclassical estimates on moments of the energy levels
\jour Acta Phys. Austriaca 
\vol 52
\yr 1980
\pages 89--105
\endref

\ref\no 9
\by E. Hewitt and K. Stromberg
\book Real and abstract analysis
\publ Springer
\publaddr New York
\yr 1965
\endref

\ref\no 10
\by E. H. Lieb and W. E. Thirring
\paper Inequalities for the moments of the eigenvalues of the Schr\"odinger
Hamiltonian and their relation to Sobolev inequalities
\inbook Studies in Mathematical Physics
\eds E. H. Lieb, B. Simon, and A. S. Wightman
\publ Princeton Univ. Press
\publaddr Princeton, NJ
\yr 1976
\pages 269--303
\endref

\ref\no 11
\by F. Mantlik and A. Schneider
\paper Note on the absolutely continuous spectrum of Sturm-Liouville operators
\jour Math. Z. 
\vol 205
\yr 1990
\pages 491--498
\endref

\ref\no 12
\by U.-W. Schmincke
\paper On Schr\"odinger{\rm '}s factorization method for 
Sturm-Liouville operators
\jour Proc. Roy. Soc. Edinburgh Sect. A
\vol 80
\yr 1978
\pages 67--84
\endref

\ref\no 13
\by V. E. Zakharov and L. D. Faddeev
\paper Kortreweg-de Vries equation{\rm :} A completely integrable Hamiltonian
system
\jour Funct. Anal. Appl.
\vol 5
\yr 1971
\pages 280--287
\endref
\endRefs
\enddocument